\documentclass[10pt,leqno]{amsart}
\usepackage{amsmath}
\usepackage{amsfonts}

\newtheorem{teor}{Theorem}[section]
\newtheorem{prop}[teor]{Proposition}
\newtheorem{coro}[teor]{Corollary}
\newtheorem{lema}[teor]{Lemma}
\newtheorem{defi}[teor]{Definition}

\theoremstyle{plain}

\def\N{\mathbb{N}}
\def\Z{\mathbb{Z}}
\def\R{\mathbb{R}}
\def\C{\mathbb{C}}

\def\todasn{\{-1,0,1\}^n}
\def\todask{\{-1,0,1\}^k}
\def\bl{\bar\lambda}
\def\bm{\bar\mu}
\def\x{\bar x}

\def\M{{\mathcal M}}
\def\D{{\mathcal D}}

\def\be{\begin{equation}}
\def\ee{\end{equation}}

\def\d{\displaystyle}

\begin{document}

\title
[Tangent, Motzkin and Catalan numbers]{Tangent and Bernoulli numbers related
to Motzkin and Catalan numbers by means of numerical triangles}
\author[J. L. Arregui]{Jos\'e Luis Arregui}
\address{Departamento de Matem\'aticas, Facultad de Ciencias, 
Universidad de Zaragoza, 50009 Zaragoza, Spain}
\email{arregui@posta.unizar.es}
\thanks{The author has been partially supported by
Proyecto  D.G.E.S. PB98-0146}



\begin{abstract} It is shown that Bernoulli numbers and tangent numbers (the derivatives of the tangent funcion at $0$) can
be obtained by means of easily defined triangles of numbers in several ways, some of them very similar to the Catalan triangle
and a Motzkin-like triangle. Our starting point in order to show this is a new expression of
$\zeta (2n)$ involving Motzkin paths.    
\end{abstract}

\maketitle

\section{Introduction. Numerical triangles generated by matrices.}

\noindent Recall that the  sequence of {\sl tangent numbers\/} is defined as 
$(\tan^{(2n-1)}(0))$, with
$$
\tan z = \sum_{n=1}^\infty \frac{\tan^{(2n-1)}(0)}{(2n-1)!}z^{2n-1} \quad (|z|<\pi/2)\,.
$$

On the other hand, the sequence of {\sl Bernoulli numbers\/} $(B_n)$ is defined by
$$
\frac z {e^z-1} = \sum_{n=0}^\infty B_n \,\frac {z^n} {n!} \quad (|z|<2\pi)\,.
$$

The first values of Bernoulli numbers are 
$$
B_0 = 1,\ B_1 = - \frac 1 2,\ B_2 = \frac 1 6,\ B_3 = 0, B_4 = -\frac 1 {30},\ B_5 =0,\ B_6 = \frac 1 {42}, 
$$
and in general $B_{2n+1} =0$ for all $n \in \N$, while the signs of $B(2n)$ alternate.

From the relation between the exponential,  sine and cosine functions, it results that 
$$
\tan^{(2n-1)}(0) = |B_{2n}| \frac {4^n\,(4^n-1)} {2n}\,.
$$

Bernoulli numbers satisfy the recurrence relation
$$
B_n = -\frac 1 {n+1}\,\sum_{k=0}^{n-1} \binom{n+1}k B_k\,.
$$

These numbers are also related to {\sl Euler numbers\/} $(E_n)$, given by $E_{2n+1}=0$ and $E_{2n}=(-1)^n E_n^*$, where
$E_n^*$ are the {\sl secant numbers\/} verifying
$$
\sec z = \frac 1 {\cos z} = \sum_{n=0}^\infty \frac {E^*_n}{(2n)!}\,z^{2n}\quad (|z|<\pi/2).
$$
Actually, if we consider $(\delta_n)$ the sequence such that
$$
\frac 1 2 \,(\sec z + \tan z) = \sum_{n=1}^\infty \delta_n z^{n-1} \quad (|z|<\pi/2)
$$
it turns out that
\begin{equation}\label{Euler}
\sum_{k=0}^\infty \frac {(-1)^{nk}}{(2k+1)^n} = \delta_n \frac{\pi^n}{2^n} \quad (n\in\N)\,
\end{equation}
a result that goes back to Euler. For even values it follows the relation between Bernoulli numbers and the Riemann
zeta function $\zeta$,
$$
\zeta(2n) = \sum_{k=1}^\infty \frac 1 {k^{2n}}= |B_{2n}|\,\frac {(2\pi)^{2n}}{2 (2n)!}\,.
$$

Standard references for all of this are \cite{GKP} and \cite{Hardy}.

A much more recent way to relate Euler and Bernoulli numbers is their combinatorial interpretation, due to R. C. Entringer,
in terms of alternating permutations. See the last section for the details. 

\bigskip

In \cite{Calabi} E. Calabi, F. Beukers and J. A. C. Kolk used a remarkable change of variables to reduce the series in (1) to
the volume of a polytope in $\R^n$, showing how it yields an elementary proof of (1). Very recently,  N. D. Elkies  has shown
in \cite{Elkies} 
how the same idea (which seems to be commonly attributed to Calabi) gives another proof of Entringer's theorem. In this paper
we follow Calabi's argument in a  different way, and although it applies only for even values of $n$ (thus
focusing on Bernoulli but not on Euler numbers) it allows us to find a nice connection between $\zeta(2n)$ and Motzkin
paths and numbers. In view of  Entringer's result, a similar connection involves Catalan numbers as well.

\bigskip

Let us first introduce one definition in order to provide an appropiate setting for our purpose: 

Let $(A^{(n)})_{n\ge 1}$ be a sequence of matrices (with numerical entries), each one with $n$ rows and $n+1$ columns. Let $t_{0,0}
=1$, $\bar t_1 = (t_{1,0}, t_{1,1}) = A^{(1)}$ and, for each $n>1$, 
$$\bar t_n = (t_{n,0}, t_{n,1}, t_{n,2},\dots,t_{n,n}) = t_{n-1} A^{(n)} = A^{(1)} A^{(2)} \cdots A^{(n)}.$$
We obtain a numerical triangle
$$
T \equiv (t_{n,m})_{0\le m\le n} \equiv \begin{array}{ccccc}
t_{0,0}  &   &   &   &   \\
t_{1,0}  & t_{1,1} &   &   &   \\
t_{2,0}  & t_{2,1} & t_{2,2} &   &   \\
t_{3,0}  & t_{3,1} & t_{3,2} & t_{3,3} &   \\
t_{4,0} & t_{4,1} & t_{4,2} & t_{4,3} & t_{4,4} \\
\hdotsfor{5}
\end{array}
$$
what we define as the {\sl triangle generated by $(A^{(n)})$}.

\medskip

In particular, if $A = (a_{ij})_{1\le i,j}$ is an infinite matrix and  $A^{(n)}$ is the submatrix formed by the first
$n$ rows and $n+1$ columns, the
triangle obtained in this way is the {\sl triangle generated by $A$}. For instance, $\bar t_3 = (t_{3,0}, t_{3,1}, t_{3,2},
t_{3,3}) = A^{(1)} A^{(2)} A^{(3)}$, that is
$$
\begin{pmatrix} a_{1,1} & a_{1,2} \end{pmatrix}
\,\begin{pmatrix} a_{1,1} & a_{1,2} & a_{1,3}  \\ a_{2,1} & a_{2,2} & a_{2,3} \end{pmatrix} \, 
\begin{pmatrix} a_{1,1} & a_{1,2} & a_{1,3} & a_{1,4} \\ a_{2,1} & a_{2,2} & a_{2,3} & a_{2,4} \\ a_{3,1} & a_{3,2} & a_{3,3} &
a_{3,4}
\end{pmatrix}\,.
$$

\bigskip

{\noindent\bf Example 1}. If $a_{ij} = 1$ for all $i,j \in\N$ the the triangle generated by $A = (a_{ij})$ is just  
$$
t_{n,m} = n! \qquad \text{for all $0\le m\le n$. }
$$

\bigskip

{\noindent\bf Example 2}. Now let $a_{ij} = \begin{cases} 2 &\text{ if} \ $i=j=1$, \\ 1 &\text{ if} \ i=j-1,\\ 0
&\text{otherwise},\end{cases}$
$$
\text{and } A = (a_{ij}) = \begin{pmatrix}
2 & 1 & 0 & 0 & \cdots  \\
0 & 0 & 1 & 0 & \cdots  \\
0 & 0 & 0 & 1 & \cdots  \\
\hdotsfor{5}\\ \end{pmatrix}.
$$
The triangle generated by $A$ is given by $t_{n,m} = 2^{n-m}$, that is
$$
\begin{array}{ccccc}
1  &   &   &   &   \\
2  & 1 &   &   &   \\
4  & 2 & 1 &   &   \\
8  & 4 & 2 & 1 &   \\
16 & 8 & 4 & 2 & 1 \\
\hdotsfor{5}
\end{array}
$$

\bigskip

{\noindent\bf Example 3}. For $a_{ij} = \begin{cases} 1 &\text{ if} \ i\le j\le i+1,\\ 0 &\text{otherwise}\end{cases}$ \quad we
have that
$$
A = (a_{ij}) = \begin{pmatrix}
1 & 1 & 0 & 0 & 0 & \cdots  \\
0 & 1 & 1 & 0 & 0 & \cdots  \\
0 & 0 & 1 & 1 & 0 & \cdots  \\
0 & 0 & 0 & 1 & 1 & \cdots  \\
\hdotsfor{6}\\ \end{pmatrix}
$$
generates the {\sl Pascal triangle}
$$
\begin{array}{cccccc}
1 &   &   &   &   &   \\
1 & 1 &   &   &   &   \\
1 & 2 & 1 &   &   &   \\
1 & 3 & 3 & 1 &   &   \\
1 & 4 & 6 & 4 & 1  &  \\
1 & 5 & 10 & 10 & 5  & 1 \\
\hdotsfor{6}
\end{array}
$$
what means that $\d t_{n,m}= \binom n m$ for all $n$ and $m$. Each number plus the one next to its left in the triangle gives
the number below. 

\bigskip

{\noindent\bf Example 4}. Take the matrix in the previous example (say the {\sl Pascal matrix\/}) and put 1 in all the entries
for the lower dominant as well, that is
$$
A = \begin{pmatrix}
1 & 1 & 0 & 0 & 0 & 0 & \cdots  \\
1 & 1 & 1 & 0 & 0 & 0 & \cdots  \\
0 & 1 & 1 & 1 & 0 & 0 & \cdots  \\
0 & 0 & 1 & 1 & 1 & 0 & \cdots  \\
0 & 0 & 0 & 1 & 1 & 1 & \cdots  \\
\hdotsfor{7}\\ \end{pmatrix}\,,
$$
Then the triangle generated by $A$ begins
$$
\begin{array}{cccccc}
1 &   &   &   &   &   \\
1 & 1 &   &   &   &   \\
2 & 2 & 1 &   &   &   \\
4 & 5 & 3 & 1 &   &   \\
9 & 12 & 9 & 4 & 1  &  \\
21 & 30 & 25 & 14 & 5  & 1 \\
\hdotsfor{6}
\end{array}
$$
This is the so-called {\sl Motzkin triangle\/}. Each number in the triangle is the sum of the one above and its (one or two)
contiguous ones.

The first column in the triangle gives the sequence of {\sl Motzkin numbers\/} $M_n$. 

\bigskip

{\noindent\bf Example 5}. If we fill with 1 all below the diagonal in the Pascal matrix  we get
$$
A = \begin{pmatrix}
1 & 1 & 0 & 0 & 0 & \cdots  \\
1 & 1 & 1 & 0 & 0 & \cdots  \\
1 & 1 & 1 & 1 & 0 & \cdots  \\
1 & 1 & 1 & 1 & 1 & \cdots  \\
\hdotsfor{6}\\ \end{pmatrix}\,,
$$
which generates
$$
\begin{array}{cccccc}
1 &   &   &   &   &   \\
1 & 1 &   &   &   &   \\
2 & 2 & 1 &   &   &   \\
5 & 5 & 3 & 1 &   &   \\
14 & 14 & 9 & 4 & 1  &  \\
42 & 42 & 28 & 14 & 5  & 1 \\
\hdotsfor{6}
\end{array}
$$
that is known as {\sl Catalan triangle\/}. Note that (if we set $t_{n,-1} = 0$ for all $n$)
$$
t_{n+1,m} = \sum_{k=m-1}^n t_{n,m}\,.
$$
The first (second) column in the triangle gives the sequence of {\sl Catalan numbers\/} $C_n$. It is well known that 
$$
t_{n,m} = \binom{2n-m} n\,\frac{m+1}{n+1} \quad (0\le m\le n)
$$
and thus $\d C_n = \binom{2n} n\,\frac 1 {n+1} = \frac {(2n)!}{n!\,(n+1)!}$.

\bigskip

\medskip

Both Motzkin and (much better known) Catalan numbers have a nice variety of realizations as the number of solutions to 
some combinatorial problem on $n$ objects (see \cite{Shap} and \cite{Gardner}). In this article we will briefly
explore  one of them,  obtaining a new one that links these numbers to permutation groups and then to Bernoulli numbers.

\bigskip

For instance, we will obtain the following theorem:

\begin{teor}\label{BerCat} Let $A = (a_{ij})$ the matrix given by

$$
a_{ij}= \begin{cases} j(j+1) &\text{ if }i \ge j-1,\\ 0 &\text{ if }i<j-1,\end{cases} 
$$
i.e.
$$
A = \begin{pmatrix}
2 & 6 & 0 & 0 & 0 & \cdots  \\
2 & 6 & 12 & 0 & 0 & \cdots  \\
2 & 6 & 12 & 20 & 0 & \cdots  \\
2 & 6 & 12 & 20 & 30 & \cdots  \\
\hdotsfor{6}\\ \end{pmatrix}\,.
$$
The first column of the triangle generated by $A$ is then the sequence of tangent numbers $\big(\tan^{(2n-1)}(0)\big)$
\end{teor}

Note that $A$ is  the product of the matrix in example 5 and (on the right) the diagonal infinite
matrix
$\big(j(j+1)
\delta_{ij}\big)$.

The triangle begins
$$
\begin{array}{cccccc}
1 &   &   &   &   \\
2 & 6 &   &   &   \\
16 & 48 & 72 &   &    \\
272 & 816 & 1440 & 1440 &   \\
7936 & 23808 & 44352 & 57600 & 43200  \\
\hdotsfor{5}
\end{array}
$$
where  $\d t_{n+1,m} = m(m+1)\sum_{k=m-1}^n t_{n,k}$.

We find it convenient to start with Motzkin paths and numbers, then showing their relation with Calabi's idea.

\bigskip\bigskip\bigskip\bigskip\bigskip

\section{Motzkin paths}

\noindent A {\sl Motzkin path\/} of $n$-th order (or $n$ steps) is a finite sequence
$$
\bar\lambda = (\lambda_1,\dots,\lambda_n) \in \todasn \text{ such that } \sum_{j=1}^n \lambda_j =0 \text { and } \sum_{j=1}^l
\lambda_j \ge 0 \text { if }l<n.
$$

The name {\sl Motzkin\/} is due to the fact that the number of Motzkin paths of $n$ steps is just the number $M_n$ in example
4, so this is one of the realization of Motzkin numbers (surely the most common by now), and the term {\sl
path\/} comes after the usual visualitazion, in an obvious 1-1 correspondence, of any such sequence with a path that joins $(0,0)$
with
$(n,0)$ in $n$ steps --in the $\Z \times
\Z$ lattice--, each one by summing the vector
$(1,\lambda_j)$, with the condition that no intermediate point lies below the horizontal axis.

The {\sl empty sequence\/} is considered as the only Motzkin path of $0$ steps, to give $M_0=1$. 

The only Motzkin path  of 1 step is $(0)$, and with 2 steps we have $(0,0)$ and $(1,-1)$. 

The four Motzkin paths of 3 steps are
\bigskip
$$
\begin{array}{cccccc}
& (0,0,0)\ \ \rightarrow\,\rightarrow\,\rightarrow & & &
\nearrow\,\searrow{{}\atop\longrightarrow}\ \ (1,-1,0) &
\\ & & & & & \\ &
(0,1,-1)\ \ {{}\atop\longrightarrow}\nearrow\,\searrow & & &
\nearrow{\longrightarrow\atop{}}\searrow\,\ (1,0,-1) & 
\end{array}
$$
\bigskip\bigskip

\noindent and a 6 step Motzkin path is $(1,0,1,-1,0,-1)$, visualized as
\bigskip
$$
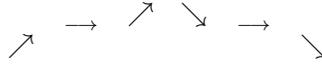
\begin{array}{cccccc}
& {{}\atop\longrightarrow} & \nearrow &\searrow& {{}\atop\longrightarrow} & \\
\nearrow& & & & &\searrow
\end{array}
$$
\bigskip

We will denote $\M_n$ the set of all Motzkin paths of $n$ steps, being then $M_n$  its cardinality. $\M = \cup_{n=1}^\infty$
$\M_n$ is the set of all Motzkin paths of any positive order. Being so easy to draw Motzkin paths, the reader can check that
 the
first Motzkin numbers are 1,1,2,4,9,21 (beginning by $n=0$).

The usual way to obtain Motzkin numbers is as in example 4. Another recursive formula is
$$
M_{n+1}= M_n + \sum_{k=0}^{n-1} M_k M_{n-k-1},
$$
and in the next section we will recover the formula that best relates them to Catalan numbers. These and further results can be
seen in \cite{Aig} and \cite{Barc}.
  
\bigskip

To fix notations, let $(\bar\lambda,\bm)=(\lambda_1,\dots,\lambda_{n_1},\mu_1,\dots,\mu_{n_2})$ if $\bar\lambda =
(\lambda_1,\dots,\lambda_{n_1})$ and $\bm=(\mu_1,\dots,\mu_{n_2})$. Here $\bm$ could be the empty sequence, and 
$(\bl,\bm) = \bl$.

The following proposition is then pictorially evident; it says that by removing flat steps, or by ``flattening cusps" in a
Motzkin path, we get another Motzkin path.

\begin{prop}\label{quitaypon} For any $n\in\N$
\item{(\text i)} $(\bl,0,\bm)\in\M_{n+1}$ if and only if $(\bl,\bm)\in\M_n$.
\item{(ii)} $(\bl,1,-1,\bm)\in\M_{n+1}$ if and only if $(\bl,0,\bm)\in\M_n$.\hfill\qed
\end{prop}

\bigskip

For each $k\in\N$, let $-1_k = (\underset{\text{$k$ times}}{-1,\dots,-1})$, so  $(\bl,-1_k) = (\bl,\underset{\text{$k$
times}}{-1,\dots,-1})$. We allow $k=0$ by defining $-1_0$ as the empty sequence. In the same way, we define $1_k$ and $0_k$
for any $k\in\N\cup\{0\}$.

\medskip

Now we define, for any $\bm\in\todask$, $\M_{n;\bm} = \{(\bl,\bm)\in \M_n\}$, the set of Motzkin paths of $n$ steps that
``end in $\bm$", and  $M_{n;\bm}$ is its cardinality.

\medskip

Note that $\M_{n;1}$ is empty, and then $\M_n = \M_{n;0} \cup \M_{n;-1}$. Moreover we have the following:

\begin{prop}\label{descom} For every $n\in\N$ 
$$
\M_n = \M_{n;0} \cup \M_{n;1,-1} \cup \M_{n;0,-1} \cup \M_{n;1,-1_2} \cup \M_{n;0,-1_2} \cup \M_{n;1,-1_3} \cup \dots
$$
(disjoint union), where the nonempty sets are exactly the  first $n$ sets in the list.
\end{prop}
\begin{proof} Note that $\M_{n;-1_k} = \M_{n;1,-1_k} \cup \M_{n;0,-1_k} \cup \M_{n;1,-1_{k+1}} \cup \dots$ is empty if $k>2n$,
since for any $(\bl,-1_k)\in \M_n$ we have $0 = \lambda_1 +\cdots+ \lambda_{n-k} - k \le n-2k$.

If $n=2k$ then $(1_k,-1_k) \in \M_{n;1,-1_k}$, while the same argument as before shows that $\M_{n;0,-1_k}$ is empty.

If $n=2k+1$ then $(1_k,0,-1_k)\in \M_{n;0,-1_k}$ and $(0,1_k,-1_k)\in\M_{n;1,-1_k}$.

Finally,  for $n>2k$  $(0_{n-2k},1_k,-1_k) \in\M_{n;1,-1_k}$ and $(0_{n-2k-1},1_k,0,-1_k)
\in\M_{n;0,-1_k}$.
\end{proof}

\bigskip

This decomposition of $\M_n$ allows us to prove the following ``general" theorem, in the context of triangles generated
by matrices.

\begin{teor}\label{Mabs} Let $f\colon\M\to\R$ a function verifying the following recursion conditions (I) and (II), in terms of
two matrices
$(b_{n,k})_{1\le n,\,0\le k}$ and $(c_{n,k})_{1\le n,k}$ with $b_{n,0}\ne 0$ for all $n$:
\item{(I)} If $(\bl,-1_k)\in \M_n$ then $f(\bl,0,-1_k) = b_{n,k} f(\bl,-1_k)$.
\item{(II)} If $(\bl,0,-1_{k-1})\in \M_n$ then $f(\bl,1,-1_k) = c_{n,k} f(\bl,0,-1_{k-1})$.

Then the sequence $\big( b_{n,0} \sum_{\bl\in\M_n} f(\bl)\big)_{n\ge 1}$ is the first column of the numerical triangle $f(0)
T$, where $T$ is the triangle generated by  the sequence of matrices $(A^{(n)})$, $A^{(n)}$ being the $n\times (n+1)$ matrix
$$
A^{(n)} = \begin{pmatrix}
b_{n,0} & c_{n,1} & 0 & 0 & 0 & 0 & \cdots  \\
b_{n,0} & 0 & b_{n,1} & 0 & 0 & 0 & \cdots  \\
b_{n,0} & 0 & b_{n,1} & c_{n,2} & 0 & 0 &\cdots  \\
b_{n,0} & 0 & b_{n,1} & 0 & b_{n,2} & 0 &\cdots  \\
b_{n,0} & 0 & b_{n,1} & 0 & b_{n,2} & c_{n,3} &\cdots  \\
\hdotsfor{7}\\
b_{n,0} & 0 & b_{n,1} & 0 & b_{n,2} & 0 &\cdots
 \end{pmatrix}\,.
$$
\end{teor}

\begin{proof} Let $(x_{n,1}, x_{n,2}, \dots , x_{n,n+1})$ be the vector given by
$$
\big(\sum_{\bl\in\M_{n+1;0}} f(\bl),\sum_{\bl\in\M_{n+1;1,-1}} f(\bl),\sum_{\bl\in\M_{n+1;0,-1}} f(\bl),\sum_{\bl\in\M_{n+1;1,-1_2}}
f(\bl),\quad\dots 
\big)
$$
i.e.
\begin{align*}
x_{n,2k+1}\,&= \sum_{\bl\in\M_{n+1;0,-1_k}} f(\bl) \qquad (0\le 2k \le n),\\
x_{n,2k}\,&= \sum_{\bl\in\M_{n+1;1,-1_k}} f(\bl) \qquad (2\le 2k \le n+1).
\end{align*}
We claim that this vector equals $f(0) A^{(1)}A^{(2)}\cdots A^{(n)}$, as we see by induction on $n$: for $n=1$ it is $(f(0,0),
f(1,-1)) = (b_{1,0}f(0), c_{1,1}f(0)) = f(0) A^{(1)}$. Assuming it true for $n-1$,  note first that
\begin{align*}
x_{n,1} &= \sum_{\bl\in\M_{n+1;0}} f(\bl) = \sum_{(\bl,0)\in\M_{n+1}} f(\bl,0) = \sum_{\bl\in\M_n}
f(\bl,0)\\ &= \sum_{\bl\in\M_n} b_{n,0}f(\bl) = b_{n,0}\sum_{\bl\in\M_n}f(\bl) = b_{n,0} \sum_{j=1}^n x_{n-1,j}\\
&= f(0) \,A^{(1)} \cdots A^{(n-1)} {\begin{pmatrix} b_{n,0}\\\cdots\\b_{n,0}\end{pmatrix}}.
\end{align*}
\noindent For $1\le k \le n/2$ we have
\begin{align*}
x_{n,2k+1} &= \sum_{\bl\in\M_{n+1;0,-1_k}} f(\bl) =  \sum_{(\bl,-1_k)\in\M_n}
f(\bl,0,-1_k)\\
&= \sum_{(\bl,-1_k)\in\M_n} b_{n,k}f(\bl,-1_k) = b_{n,k}\sum_{\bl\in\M_{n;-1_k}}f(\bl)\\
&= b_{n,k} \big(\sum_{\bl\in\M_{n;1,-1_k}}\hskip-7pt f(\bl) + \sum_{\bl\in\M_{n;0,-1_k}}\hskip-7pt f(\bl) +
\sum_{\bl\in\M_{n;1,-1_{k+1}}}\hskip-7pt f(\bl) +\cdots \big)\\ &= b_{n,k} \sum_{j=2k}^n x_{n-1,j} = f(0) A^{(1)} \cdots
A^{(n-1)} {\begin{pmatrix} 0\\\cdots\\0\\b_{n,k}\\\cdots\\b_{n,k}\end{pmatrix}}\hskip-7pt\begin{array}{l}
   \\
   \\
   \\
{}^{{}^{\leftarrow 2k}}   \\
  \\
\end{array},
\end{align*}
and finally, for $1 \le k \le \dfrac {n+1}2$,
\begin{align*}
x_{n,2k} &= \sum_{\bl\in\M_{n+1;1,-1_k}} f(\bl) =  \sum_{(\bl,0,-1_{k-1})\in\M_n}
f(\bl,1,-1_k)\\
&= \sum_{(\bl,0,-1_{k-1})\in\M_n} c_{n,k}f(\bl,0,-1_{k-1}) =c_{n,k}\sum_{\bl\in\M_{n;0,-1_{k-1}}}f(\bl)\\
&= c_{n,k}\, x_{n-1,2k-1} = f(0) A^{(1)} \cdots A^{(n-1)} {\begin{pmatrix}
0\\\cdots\\0\\c_{n,k}\\0\\\cdots\\0\end{pmatrix}}\hskip-7pt\begin{array}{l}
   \\
   \\
   \\
{}^{{}^{\leftarrow 2k-1}}   \\
  \\
  \\
\end{array}\,,\\
\end{align*}
so our claim is true.
\end{proof}
Our first application of this theorem is  counting Motzkin paths.

\begin{coro}\label{Motzkin} The sequence of Motzkin numbers $(M_n)$ is the first column of the triangle generated by $A = (a_{ij})$
with
$$
a_{ij}= \begin{cases} 1 & \text{ if } i=j-1, \\ 1 &\text{ if } j \text{ is odd and } i>j-1,\\ 0 &otherwise\,.\end{cases} 
$$
\end{coro}
\begin{proof} Apply Theorem \ref{Mabs} to $f = 1$. \end{proof}

\medskip

\noindent Note that 
$$
A = \begin{pmatrix}
1 & 1 & 0 & 0 & 0 & 0 &  \cdots  \\
1 & 0 & 1 & 0 & 0 & 0 &  \cdots  \\
1 & 0 & 1 & 1 & 0 & 0 &  \cdots  \\
1 & 0 & 1 & 0 & 1 & 0 &  \cdots  \\
1 & 0 & 1 & 0 & 1 & 1 &  \cdots  \\
1 & 0 & 1 & 0 & 1 & 0 &  \cdots  \\
\hdotsfor{7} \end{pmatrix}\,,
$$
and the triangle (which is not the Motzkin triangle of example 4) begins
$$
\begin{array}{cccccc}
1 &   &   &   &   &   \\
1 & 1 &   &   &   &   \\
2 & 1 & 1 &   &   &   \\
4 & 2 & 2 & 1 &   &   \\
9 & 4 & 5 & 2 & 1  &  \\
21 & 9 & 12 & 5 & 3  & 1 \\
\hdotsfor{6}
\end{array}
$$
Once the $(n-1)$-th row is given, we get the number $t_{n,m}$ as in Catalan triangle for odd $m$, but $t_{n,m} = t_{n-1,m-1}$
if
$m$ is even.

Following the proof of Theorem \ref{Mabs}, note that the $n$-th row is just
$$
(M_{n+1;0}, M_{n+1;1,-1}, M_{n+1;0,-1}, M_{n+1;1,-1_2},\cdots)\,.
$$
For instance, the third row is $(4,2,2,1)$, the cardinalities of
\begin{align*}
\M_{4;0}&=\{(0,0,0,0),(1,-1,0,0),(1,0,-1,0),(0,1,-1,0)\}\,,\\
\M_{4;1,-1}&=\{(0,0,1,-1),(1,-1,1,-1)\}\,,\\ 
\M_{4;0,-1}&=\{(1,0,0,-1),(0,1,0,-1)\}\,,\\
\M_{4;1,-1_2}&=\{(1,1,-1,-1)\}\,,
\end{align*}
which sum up $9 = M_4 =M_{5;0} = t_{4,0}$.

\bigskip\bigskip\bigskip\bigskip\bigskip

\section{Dyck paths and Catalan numbers.}

\noindent There is antoher easy decomposition of Motzkin paths that relates Motzkin and Catalan numbers, by counting how many
$\lambda_j = 0$ there are in each $\bl$. For each
$n\in\N$ and
$k=0,1,\dots,n$, let 
$$
\D_{n,k} = \{\bl \in\M;\ \bl \text{ has } k \text { 0's}\}
$$
and let $D_{n,k}$ its cardinality. $D_{n,k}$ is then the number of Motzkin paths of $n$ steps, $k$ of which are horizontal.

\bigskip

Since any Motzkin path has the same number of 1's and $-$1's, the number of horizontal steps is of the same parity as $n$.
Hence
$D_{n,k}=0$ if $n$ is even and $k$ is odd or vice-versa.

\bigskip

It is also obvious that $D_{n,n} =1$ for each $n$.

\bigskip

$\D_{n,0}$ is the set of the so-called {\sl Dyck paths\/}, i.e. Motzkin paths with no horizontal steps. It is empty for
any
 odd $n$, and then 
$$
\D_{0}=\cup_{n=1}^\infty \D_{2n,0}
$$
is the set of Dyck paths of any order.
It is well known that 
$$
D_{2n,0}=C_n
$$
for all $n$, where $C_n$ is the Catalan number of example 5. See for instance \cite{Dyck}. We will find out a simple proof of
this.

\bigskip

With the same notations as in the previous section, we clearly have
$$
\D_{2n,0} = \D_{2n,0;1,-1} \cup \D_{2n,0;1,-1_2} \cup \dots \cup \D_{2n,0;1,-1_n}
$$
(disjoint union). Denote the cardinality of each set by $D_{2n,0;1,-1_k} >0$.

\bigskip

Moreover, note that
$$
(\bl,-1_k)\in \D_{2n,0;-1_k} \text{ if and only if } (\bl, 1, -1_{k+1}) \in \D_{2n+2,0;1,-1_{k+1}}\,,
$$
and
$$
\D_{2n,0;-1_k} = \D_{2n,0;1,-1_k}\cup\D_{2n,0;1,-1_{k+1}}\cup\dots\cup\D_{2n,0;1,-1_n}\,.
$$

Similarly to Theorem \ref{Mabs}, we can formulate a ``general" theorem:

\begin{teor}\label{Cabs} Let $f\colon\D_0\to\R$ a function verifying the following recursion condition, in terms of 
a matrix
$(a_{n,k})_{1\le n,\,0\le k}$  with $a_{n,0}\ne 0$ for all $n$:
\item{} If $(\bl,-1_k)\in \D_{2n,0}$ then $f(\bl,1,-1_{k+1}) = a_{n,k} f(\bl,-1_k)$.

Then the sequence $\d\big( a_{n,0} \sum_{\bl\in\D_{2n,0}} f(\bl)\big)_{n\ge 1}$ is the first column of the numerical triangle
$f(1,-1) T$, where $T$ is the triangle generated by  the sequence of matrices $(A^{(n)})$, $A^{(n)}$ being the $n\times (n+1)$
matrix
$$
A^{(n)} = \begin{pmatrix}
a_{n,0} & a_{n,1} & 0 & 0 & 0 & \cdots & 0 \\
a_{n,0} & a_{n,1} & a_{n,2} & 0 & 0 & \cdots & 0 \\
a_{n,0} & a_{n,1} & a_{n,2} & a_{n,3} & 0 & \cdots & 0  \\
a_{n,0} & a_{n,1} & a_{n,2} & a_{n,3} & a_{n,4} &\cdots & 0  \\
\hdotsfor{7}\\
a_{n,0} & a_{n,1} & a_{n,2} & a_{n,3} & a_{n,4} &\cdots &a_{n,n}
 \end{pmatrix}\,.
$$
\end{teor}

\begin{proof} The proof is similar to that of Theorem \ref{Mabs}. We have to prove that
$$
\big( \sum_{\bl\in\D_{2n+2,0;1,-1}} \hskip-7pt f(\bl), \sum_{\bl\in\D_{2n+2,0;1,-1_2}} \hskip-7pt f(\bl), \ \dots,
\sum_{\bl\in\D_{2n+2,0;1,-1_{n+1}}} \hskip-7pt f(\bl)\ \big)
$$ 
equals $f(1,-1) A^{(1)} A^{(2)} \cdots A^{(n)}$,  
the first component being $a_{n,0} \d\sum_{\bl\in\D_{2n,0}} f(\bl)$. \end{proof}

Take in particular $f=1$. Then $a_{n,k}=1$ for all $0\le k \le n$, and $A$ is just the matrix of example 5. Theorem
\ref{Cabs} says in this case that
$$
C_n = \sum_{\bl\in\D_{2n,0}} 1 = D_{2n,0}\,.
$$

As in ``our" Motzkin triangle, we also get that the numbers of the $n$-th row in the Catalan triangle are the number of Dyck
paths of $2n$ steps finishing in $-1_k$, for $k=0,\dots,n$. For instance, the third row is 
$$(5,5,3,1) = (D_{6,0}, D_{6,0;-1},
D_{6,0;-1_2}, D_{6,0;-1_3})\,,$$ the cardinalities of (listed backwards)
\begin{align*}
\D_{6,0;-1_3} &= \{(1,1,1,-1,-1,-,1)\}\,,\\
\D_{6,0;-1_2} &= \D_{6,0;-1_3} \cup \{(1,-1,1,1,-1,-1),(1,1,-1,1,-1,-1)\}\,,\\
\D_{6,0;-1} &= \D_{6,0;-1_2} \cup \{(1,-1,1,-1,1,-1),(1,1,-1,-1,1,-1)\}\,,\\
\D_{6,0} &= \D_{6,0;-1}\,.
\end{align*}

\bigskip

Our results concerning Bernoulli numbers are applications of Theorems \ref{Mabs} and \ref{Cabs}. Before seing that, we show
another application of Theorem \ref{Mabs}:

\begin{teor}\label{fceros} Let 
$$
A(x) = \begin{pmatrix}
x & 1/x & 0 & 0 & 0 & 0 &  \cdots  \\
x & 0 & x & 0 & 0 & 0 &  \cdots  \\
x & 0 & x & 1/x & 0 & 0 &  \cdots  \\
x & 0 & x & 0 & x & 0 &  \cdots  \\
x & 0 & x & 0 & x & 1/x &  \cdots  \\
x & 0 & x & 0 & x & 0 &  \cdots  \\
\hdotsfor{7} \end{pmatrix}\,,
$$
and let $T(x)$ the triangle generated by $A(x)$.  Then its first column is given by the polynomials
$$
\sum_{k=0}^n D_{n,k}\, x^k\,.
$$
\end{teor}
\begin{proof}  Let $z(\bl)$ the function that counts the zeroes in each $\bl\in\M$, and let $f(\bl) = x^{z(\bl)}$. Note that
\begin{align*}
f(\bl,1,-1_k) &= \frac 1 x f(\bl, 0, -1_{k-1}) \text{ and }\\
f(\bl,0,-1_k) &= x f(\bl, -1_k)\,,
\end{align*}
so Teorem \ref{Mabs} applies, and the matrices $A^{(n)}$ for each $x$ are the corresponding submatrices of $A(x)$. According to
this, the first element of the $n$-th row is
$$
\sum_{\bl\in\M_n} x^{z(\bl)} = \sum_{k=0}^n D_{n,k}\, x^k\,.
$$
\end{proof}

Some remarks are in order: of course, if $x=1$ it results the matrix of Corollary \ref{Motzkin}, what is just saying that $M_n
=
\sum_{k=0}^n D_{n,k}$, obvious since $(\D_{n,k})_{0\le k\le n}$ give a decomposition of $\M_n$.

\bigskip

On the other hand, note that from each Dyck path of $\D_{2m-2k,0}$ we get, by interlacing $2k$ 0's in all the
$\binom{2m}{2k}$ possible ways, this number of different paths in $\D_{2m,2k}$, and any element in this set is obtained
in this way from a unique path in
$\D_{2m-2k,0}$. It follows that
$$
D_{2m,2k}=\binom{2m}{2k} \, C_{m-k} = \binom{2m}{2(m-k)} \, C_{m-k}\,.
$$
Analogously
$$
D_{2m+1,2k+1}=\binom{2m+1}{2k+1} \, C_{m-k} = \binom{2m+1}{2(m-k)} \, C_{m-k}\,.
$$

These two facts give the most usual relation between Motzkin and Catalan numbers,
$$
M_n = \sum_{n\ge 0} \binom n{2k} \, C_k\,.
$$ 

\bigskip\bigskip\bigskip\bigskip\bigskip

\section{Calabi's argument and Motzkin paths.}

\noindent Let us briefly explain Calabi's argument reducing the
calculation of $\zeta(2n)$ to the volume of a polytope in $\R^n$:

For any natural number $n\ge 2$, by writing
$$
\zeta(n) = \sum_{k=1}^\infty \frac 1 {k^n} = \sum_{k=0}^\infty\frac 1 {(2k+1)^n} + \sum_{k=1}^\infty \frac 1 {(2k)^n}
$$
we get 
$$
\zeta(n) (1-2^{-n}) = \sum_{k=0}^\infty\frac 1 {(2k+1)^n}.
$$
Since $\d\frac 1 {2k+1} = \int_0^1 x^{2k}dx$, then
$$
\zeta(n) = \frac {2^n}{2^n-1} \sum_{k=0}^\infty \int_{(0,1)^n}(x_1 x_2 \cdots x_n)^{2k} dm(\x),
$$
where of course $\x = (x_1,\dots,x_n)$ and $m$ stands for the Lebesgue measure. Then
\begin{align*}
\zeta(n) &= \frac {2^n}{2^n-1} \int_{(0,1)^n} \sum_{k=0}^\infty(x_1 x_2 \cdots x_n)^{2k} dm(\x)\\ &= \frac {2^n}{2^n-1}
\int_{(0,1)^n}
\frac 1 {1-x_1^2 x_2^2 \cdots x_n^2} dm(\x)\,.
\end{align*}
Try now the change of variables
$$
(x_1,x_2,\dots,x_n) = \varphi(\bar u) =\big(\frac {\sin u_1}{\cos u_2},\frac {\sin u_2}{\cos u_3},\dots,
\frac {\sin u_{n-1}}{\cos u_n}, \frac {\sin u_n}{\cos u_1}\big)\,.
$$
Note that if $0<u,v<\pi/2$ then $\d\frac{\sin u}{\cos v} = \frac{\sin u}{\sin(\pi/2 - v)}>0$, and it is $<1$ if and only if
$u<\pi/2 -v$. Then, if we define
$$
A_n = \{\bar u >0;\ u_1 + u_2 < \pi/2, u_2 + u_3 < \pi/2, \dots, u_{n-1}+u_n < \pi/2, u_n + u_1 < \pi/2\} 
$$
we have that $\varphi(A_n) \subseteq (0,1)^n$. It is not difficult to see that  $\varphi(A_n)$ is actually $(0,1)^n$ (moreover
$\varphi$ is injective on $A_n$), and hence this change of variables gives
$$
\int_{(0,1)^n}
\frac 1 {1-x_1^2 x_2^2 \cdots x_n^2} dm(\x) = \int_{A_n} \frac {\det (J\varphi(\bar u))}{1-\tan^2 u_1  \cdots \tan^2 u_n}
dm(\bar u).
$$ 
Now, it is readily checked that
$$
\det (J\varphi(\bar u)) = 1 + (-1)^{n+1}\tan^2 u_1 \cdots \tan^2 u_n,
$$
and thus, if $n$ is even, it results
$$
\int_{(0,1)^n}
\frac 1 {1-x_1^2 x_2^2 \cdots x_n^2} dm(\x) = \int_{A_n} dm = m(A_n).
$$
Take instead
$$
E_n = \{0<\x\in\R^n; x_1 + x_2 <1, x_2 + x_3 <1, \dots, x_{n-1}+x_n <1, x_n + x_1 <1\}\,.
$$ 
Since
$m(A_n) =
\dfrac {\pi^n}{2^n} m(E_n)$, we have shown that
\begin{equation}\label{politopo}
\zeta(2n) = \frac{\pi^{2n}}{4^n-1}m(E_{2n}) 
\end{equation}
for all $n\in\N$.
Now let $x_{2j}=y_{2j}$ and $x_{2j-1} = 1-y_{2j-1}$ for $j=1,\dots,n$. This new change of variables gives
$$
\zeta(2n) = \frac{\pi^{2n}}{4^n-1}m(F_{2n})\,, 
$$
where $F_{2n} = \{ 0<\bar y <1;\ y_2<y_1,\, y_2<y_3, \dots, y_{2n}<y_{2n-1},\, y_{2n}<y_1\}$.

\bigskip

For $n=1$, we have $F_2 = \{(a,b)\in\R^2;\ 0<b<a<1\}$, so $m(F_2) = 1/2$ and $\d\zeta(2) =  {\pi^2}/6$.

\bigskip

Changing the coordinates order, we can re-write 
$$
F_{2n} = \{(\bar t,\x) \in (0,1)^n\times(0,1)^n;\ t_1<x_1,\,t_1<x_2,\dots, t_n< x_n,\, t_n<x_1\}.
$$
Let $a\wedge b = \min\{a,b\}$ for all real numbers $a$ and $b$. 
For each $\x\in (0,1)^n$ the corresponding section of $F_{2n}$ is
\begin{align*}
F_{\x} &= \{\bar t;\ (\bar t,\x) \in F_{2n}\} \\&= (0,x_1\wedge x_2) \times (0,x_2\wedge x_3)
\times\cdots\times(0,x_{n-1}\wedge x_n)\times(0,x_n\wedge x_1),
\end{align*}
so
$$
m(F_{2n}) = \int_{(0,1)^n}\hskip-5pt m(F_{\x}) dm(\x) = \int_{(0,1)^n}(x_1\wedge x_2)\cdots(x_{n-1}\wedge x_n)(x_n\wedge
x_1)dm(\x),
$$
and we have obtained the following expression of $\zeta(2n)$ for each natural $n$:
\begin{equation}\label{xi}
\zeta(2n) = \frac{\pi^{2n}}{4^n-1}\int_{(0,1)^n}\xi(\x)dm(\x)\,, 
\end{equation}
where $\xi(\x)=(x_1\wedge x_2)(x_2\wedge x_3)\cdots(x_{n-1}\wedge x_n)(x_n \wedge x_1)$.

\bigskip

For $n=2$  we have $\xi(a,b)=(a\wedge b)^2$, and 
\begin{align*}
\int_{(0,1)^2} (a\wedge b)^2 dm(a,b) &= 2 \int_{\{0<b<a<1\}} \hskip-10pt(a\wedge b)^2 dm(a,b) = 2 \int_{\{0<b<a<1\}}
\hskip-10pt b^2 dm(a,b)\\ &= 2\int_0^1\int_0^a b^2 db da = \frac 2 3 \int_0^1 a^3 da = \frac 1 6\,,
\end{align*}
so (3) says that $\d\zeta(4) = \frac{\pi^4}{90}$.

\bigskip

For $n=3$ we can also write
$$
m(D_6) = 6 \int_{\{0<x_1<x_2<x_3<1\}}\hskip-10pt\xi(\x)dm(\x),
$$ 
since any coordinate is compared to each other in $(x_1\wedge x_2) (x_2\wedge x_3) (x_3\wedge x_1)$. This makes it
easy to see that $m(F_6) = 1/15$ and thus $\zeta(6) = \pi^6/945$.

\bigskip

For $n>3$ it does not suffice to consider the set $\{0<x_1<\cdots <x_n<1\}$.

\bigskip

Let $S_n$ be the group of permutations of the set $\{1,2,\dots,n\}$, for each $n\in\N$. 

Note that, in the integral in (3), we can ignore all points in $(0,1)^n$ with two equal coordinates, since $m(H)=0$ for any
hyperplane
$H$. Besides, if we fix $j$ and integrate in $\{\x; \ x_j<x_i \text{ for all }i\ne j,\ x_i\ne x_k\text{ for all }i,k\}$ then the
result is independent from
$j$, by simmetry. Choose $j=n$. For each $\x$ in the corresponding set there exists a unique permutation $\sigma\in S_{n-1}$ such
that
$$
0<x_n<x_{\sigma(1)}<x_{\sigma(2)}<\cdots<x_{\sigma(n-1)}<1.
$$
If $E_\sigma$ is the set of such points for  any given $\sigma$, we have 
\begin{equation}\label{trozos}
\int_{(0,1)^n}\xi(\x) dm(\x) = n \sum_{\sigma\in S_{n-1}} \int_{E_\sigma}\xi(\x) dm(\x) = n \sum_{\sigma^{-1}\in S_{n-1}}
\int_{E_{\sigma^{-1}}}\xi(\x) dm(\x)\,.
\end{equation}
Let $E=E_{\text{id}} = \{\x;\ 0<x_n <x_1 <x_2 <\cdots < x_{n-1}<1\}$.

The change of variables $\x\mapsto (x_{\sigma(1)},\dots, x_{\sigma(n-1)}, x_n)$ maps $E$ onto $E_{\sigma^{-1}}$,
and then
\begin{align*}
\int_{E_{\sigma^{-1}}}\xi(\x) dm(\x) &= \int_E \xi(x_{\sigma(1)},\dots,x_{\sigma(n-1)},x_n)dm(\x)\\
&= \int_E x_n^2 (x_{\sigma(1)}\wedge x_{\sigma(2)})\cdots (x_{\sigma(n-2)}\wedge x_{\sigma(n-1)}) dm(\x)\,.
\end{align*}

\bigskip\bigskip

The last integral is easily expressed in terms of the ad hoc map
$$
\Phi\colon \cup_{n=2}^\infty S_n\longrightarrow \cup_{n=1}^\infty \{-1,0,1\}^n
$$
defined as follows: for any $n\in\N$ and $\sigma\in S_{n+1}$, $\Phi(\sigma)$ is the family $\bl = (\lambda_1,\dots,\lambda_n)\in
\{-1,0,1\}^n$ such that, in the list
$$0,\sigma(1), \sigma(2), \dots, \sigma(n), \sigma(n+1), 0
$$ 
$\lambda_j +1$ is the number of neighbours of $j$ that are greater than $j$.

\bigskip

For example, if $\sigma\in S_5$ is the given by $0,2,1,4,5,3,0$ then $\Phi(\sigma)=(1,-1,0,0)$: $\lambda_1 +1 =2$ because 2 and 4
(the neighbours of 1) are both $>1$; $\lambda_2 +1 =0$ since $0,1<2$, and $\lambda_3 +1 = \lambda_4 +1 =1$, for both 3 and 4 have 5
as a neighbour.

\bigskip

It is easy that $\sum_{j=1}^n \lambda_j =0$ for all $\bl = \Phi(\sigma)$: consider all the couples of neighbours; appart from the
first and the last (containing 0) each one adds $+1$ to $\sum_{j=1}^n \lambda_j -n$, and there are $n$ such couples.

\medskip

On the other hand it is clear that $\lambda_1 > -1$ and $\lambda_n < 1$, as in Motzkin paths.

\bigskip

\begin{defi} For each $\bl \in \cup_{n=1}^\infty \{-1,0,1\}^n$, $\nu(\bl)$ is the number of permutations $\sigma$ such that
$\Phi(\sigma) = \bl$.
\end{defi} 

Note that $\nu(\bl) >0$ if and only if $\bl$ is in the image of $\Phi$.

\bigskip

In the following $e_j$ will denote $(0,0,\dots,0,\underset{{}^j}{1},0,\dots)$, and then $\bl \pm e_j =
(\lambda_1,\lambda_2,\dots,\lambda_j\pm 1,\lambda_{j+1},\dots)$.

\bigskip

\begin{lema}\label{simple} For each $\bl\in\{-1,0,1\}^n$ it holds:
\item{(i)} $\d\nu(\bl,0) = 2\nu(\bl)\,$;
\item{(ii)} $\d\nu(\bl,-1) = 2 \sum_{\lambda_j=0} \nu(\bl -e_j) + \sum_{\lambda_j=1}\nu(\bl - e_j)\,$.
\end{lema}
\begin{proof} Let $\sigma\in S_{n+1}$ given by $0,\sigma(1),\sigma(2),\dots,\sigma(n),\sigma(n+1),0$. From $\sigma$ we obtain
$n+2$ permutations $\tau\in S_{n+2}$, by interlacing $n+2$ between two any neighbours. These $\tau$ are different each other,
and the
$(n+2)!$ permutations of $S_{n+2}$ arise in this way by taking the $(n+1)!$ permutations of $S_{n+1}$.

If $\tau$ is so-derived from $\sigma$, with $\Phi(\sigma) = \bl$ and $\Phi(\tau)=\bm$,
then or $\lambda_j = \mu_j$ for $j=1,\dots,n$ with $\mu_{n+1}=0$ (this happens in two cases, setting $n+2$ next to $n+1$) or
$\mu_j = \lambda_j+1$ for exactly one $j$, $\lambda_j = \mu_j$ for the rest and $\mu_{n+1}=-1$.

If $\lambda_j=-1$ then $\mu_j =0$ for each $\tau$ that puts $n+2$ next to $j$, and if $\lambda_j =0$ then $\mu_j =1$ only if
$\tau$ sets $n+2$ between $j$ and its least neighbour. \end{proof}

\begin{teor} \label{imagen} The image of $\Phi$ is $\M$, the set of Motzkin paths.
\end{teor} 
\begin{proof} It is easy for $n=1$: $\M_1 = \{(0)\}$, $S_2 = \{(0,1,2,0), (0,2,1,0)\}$ and $\nu(0) = 2$. Assume that it is
true for $n$.

If $\bl \in \Phi(S_{n+2})$ and $\bl =(\bm,0)$, by Lemma \ref{simple} $\nu(\bm)\ne 0$, hence $\bm\in\M_n$ and then $\bl\in
\M_{n+1}$ (by Lemma \ref{quitaypon}). If $\bl = (\bm,-1)$ Lemma \ref{simple} says that, for some $j$, $\bm - e_j$ has preimages in
$S_{n+1}$, so $\bm-e_j\in \M_n$ and this easily implies that $\bl\in \M_{n+1}$ (a simple draw may help).

\medskip

Conversely, if $\bl\in\M_{n+1}$ and $\bl=(\bm,0)$ then $\mu\in\M_n$, and thus $2\nu(\bm) = \nu(\bl)>0$, so $\bl\in
\Phi(S_{n+2})$.

If $\bl \in \M_{n+1}$ with $\bl=(\bm,-1)$, by the definition of Motzkin paths we have $\sum_{j=1}^n \mu_j =1$ and $\sum_{j=1}^k
\mu_j \ge 0$ for $1\le k$. Let $l$ the maximum such that $\sum_{j=l}^n \mu_j =1$. As $\sum_{j<l}\mu_j =0$ it can't be $\mu_l=-1$,
so we can take $\bm - e_l$; Would it be $\sum_{j=l}^k \mu_j = -1$ for some $k>l$, it should be $\sum_{j=k+1}^n \mu_j = 2$, in
contradiction with our election of $l$. It follows that $\bm - e_l \in \M_n$, and hence $\nu(\bm - e_l)>0$. Lemma \ref{simple}
gives now that $\nu(\bm,-1) = \nu(\bl) >0$, and hence $\lambda \in \Phi(S_{n+2})$.\end{proof}

\bigskip

Let's get back to formula (4), where $n\ge 3$, and the integral 
$$
\int_E x_n^2 (x_{\sigma(1)}\wedge x_{\sigma(2)})\cdots (x_{\sigma(n-2)}\wedge x_{\sigma(n-1)}) dm(\x)\,.
$$
If $\bl = \Phi(\sigma) \in \M_{n-2}$ then this integral is the same as
\begin{align*}
\int_E x_n^2\, &x_1^{\lambda_1+1}\cdots x_{n-2}^{\lambda_{n-2}+1} dm(\x)\\
&=\int_0^1\int_0^{x_{n-1}}x_{n-2}^{\lambda_{n-2}+1}\cdots\int_0^{x_2}x_1^{\lambda_1+1}\int_0^{x_1}x_n^2\,dx_n\,dx_1\,\cdots
dx_{n-2}\,dx_{n-1}\\
&=\frac 1 {6n} \,\cdot\,\frac 1 {\lambda_1 +5} \,\cdot\,\frac 1 {\lambda_1+\lambda_2 + 7} \,\cdots\, \frac 1 {\lambda_1 +
\cdots +\lambda_{n-2} + (2n-1)}\,\,.
\end{align*}

\bigskip

\begin{defi} For each $\bl\in\todasn$, 
$$\rho(\bl) = (\lambda_1 +5) (\lambda_1 + \lambda_2 +7) \cdots (\lambda_1 + \cdots +\lambda_n
+ 2n+3)\,.$$
\end{defi}

In this expression, each $j$-th factor ($j>1$) is the previous one plus 3 (if $\lambda_j=1$), plus 2 (if $\lambda_j=0$) or plus 1
(if $\lambda_j = -1$). If $\bl\in\M$ then $\sum_j\lambda_j=0$, so the last factor is $2n+3$.

\bigskip

Summarizing, we have that
$$
\int_{(0,1)^n} \xi(\x) dm(\x) = \frac 1 6 \sum_{\sigma\in S_{n-1}} \frac 1 {\rho(\Phi(\sigma))}\,
$$
and then
\begin{equation}\label{xim}
\int_{(0,1)^n} \xi(\x) dm(\x) = \frac 1 6 \sum_{\bl\in\M_{n-2}}\frac {\nu(\bl)}{\rho(\bl)}\,,
\end{equation}
which along with (2) gives, for every $n\ge 3$,
\begin{equation}\label{zetam}
\zeta(2n) = \frac {\pi^{2n}}{6(4^n-1)}\,\sum_{\bl\in\M_{n-2}}\frac {\nu(\bl)}{\rho(\bl)}\,.
\end{equation}

\medskip

Stop for a moment and compute $\zeta(6)$ and $\zeta(8)$:

\bigskip

$\M_1 = \{(0)\}$, with $\nu(0) = 2$ and $\rho(0) = 5$. Then
$$
\zeta(6) = \frac {\pi^6}{6(4^3-1)}\,\cdot \frac 2 5 = \frac {\pi^6}{945}\,.
$$

\medskip

$\M_2 = \{(0,0),(1,-1)\}$. From Lemma \ref{simple} $\nu(0,0) = 2\nu(0) = 4$ and $\nu(1,-1)=\nu(0) =2$. Besides $\rho(0,0) = 5\cdot
7$ and $\rho(1,-1) = 6\cdot 7$, whence
$$
\zeta(8) = \frac {\pi^8}{6(4^4-1)}\,\cdot \big(\frac 4 {5\cdot 7} + \frac 2 {6\cdot 7}\big)= \frac {\pi^8}{9450}\,.
$$

\bigskip

\begin{prop} \label{general} For any $\bl\in\cup_{n=1}^\infty\todasn$ we have that 
\item{(i)} $\nu(\bl,0,-1_k) = 2 (k+1) \nu(\bl, -1_k)$ for each $k\ge 0$, and
\item{(ii)} $\nu (\bl,1,-1_k) = k (k+1) \nu(\bl, -1_{k-1}) = \dfrac{k+1}2 \nu(\bl,0,-1_{k-1})$ for each $k\ge 1$.
\end{prop}
\begin{proof} The proof follows by induction on $k$, using Lemma \ref{simple}.

\noindent(i) If $k=0$ the result is just (i) in Lemma \ref{simple}. If $k>0$ and we know that it is true for $k-1$, by $(ii)$ 
in the lemma we can write
\begin{align*}
\nu(&\bl,0,-1_k) = \nu(\bl,0,-1_{k-1},-1) \\
&= 2\nu(\bl,-1,-1_{k-1}) + 2\sum_{\lambda_j=0}\nu(\bl-e_j,0,-1_{k-1})+\sum_{\lambda_j=1}\nu(\bl-e_j,0,-1_{k-1})\\
&= 2\nu(\bl,-1_k) + 2k \big(2\sum_{\lambda_j=0}\nu(\bl-e_j,-1_{k-1})+\sum_{\lambda_j=1}\nu(\bl-e_j,-1_{k-1})\big)\\
&= 2\nu(\bl,-1_k) + 2k \nu(\bl,-1_{k-1},-1) \\
&= 2(k+1)\nu(\bl,-1_k).
\end{align*} 
\bigskip
\noindent(ii) If $k=1$ then 
\begin{align*}
\nu(\bl,1,-1) &= 2\sum_{\lambda_j=0}\nu(\bl-e_j,1) + \sum_{\lambda_j=1}\nu(\bl-e_j,1)+\nu(\bl,0)\\
&= \nu(\bl, 0) = 2\nu(\bl)
\end{align*}
since $\nu(\bm,1)=0$ for any $\bm$.

If $k>1$ and the result is true for $k-1$, then
\begin{align*}
\nu(&\bl,1,-1_k) = \nu(\bl,1,-1_{k-1},-1) \\
&= 2\sum_{\lambda_j=0}\nu(\bl-e_j,1,-1_{k-1})+\sum_{\lambda_j=1}\nu(\bl-e_j,1,-1_{k-1})+\nu(\bl,0,-1_{k-1})\\
&=(k-1)k\big(2\sum_{\lambda_j=0}\nu(\bl-e_j,-1_{k-2})+\sum_{\lambda_j=1}\nu(\bl-e_j,-1_{k-2})\big)\\
&\quad+2k\nu(\bl,-1_{k-1}) \quad\text{ (by (i))}\\
&=(k-1)k\,\nu(\bl,-1_{k-2},-1) + 2k\,\nu(\bl,-1_{k-1})\\ &= (k+1)k\,\nu(\bl, -1_{k-1}).
\end{align*}

The rest follows from (i): $\nu (\bl,0,-1_{k-1}) = 2k \nu(\bl,-1_{k-1})$, so
$$
\nu(\bl,1,-1_k) = (k+1) k \frac 1{2k} \nu(\bl,0,-1_{k-1}) = \frac{k+1}2 \nu(\bl,0,-1_{k-1}).
$$
\end{proof}

\begin{coro}\label{factM}
Let $A =(a_{ij})$ the matrix given by
$$
a_{ij} = \begin{cases} 2j &\text{ if }j \text{ is odd and }i \ge j-1\,,\\
              \dfrac{(j+2)}4 &\text{ if }j \text{ is even and }i= j-1\,,\\
            0 &otherwise,\end{cases}
$$
i.e.
$$
A = \begin{pmatrix}
2 & 1 & 0 & 0 & 0 & 0 & \cdots  \\
2 & 0 & 4 & 0 & 0 & 0 &  \cdots  \\
2 & 0 & 4 & 3/2 & 0 & 0 &  \cdots   \\
2 & 0 & 4 & 0 & 6 & 0 & \cdots  \\
2 & 0 & 4 & 0 & 6 & 2 & \cdots   \\
\hdotsfor{7}
\end{pmatrix}\,.
$$
Then, if $T =(t_{n,m})$ is the triangle generated by $A$, for all $n\in\N$ 
$$
t_{n,0} = (n+1)!\,.
$$
\end{coro}

\begin{proof} Proposition \ref{general} says that we can apply Theorem \ref{Mabs} to $f = \nu$, and conditions (I) and (II)
there are satisfied with $b_{n,k}$ and $c_{n,k}$ independent from $n$. It results that the corresponding $A^{(n)}$ matrices
are the submatrices of
$A$ as stated, and since $f(0) = 2 = b_{2,0}$ we get
$$
t_{n,0} = \sum_{\bl\in\M_n} \nu(\bl)\,,
$$
which is just the cardinality of $S_{n+1}$, that is $(n+1)!$.
\end{proof}

\bigskip

\begin{prop}\label{ro} (i) If $k\ge0$ and $(\bl,-1_k)\in \M_n$ then
\begin{align*}
\rho(\bl,0,-1_k) &=\rho(\bl) (2n+5-k) (2n+6-k) \cdots (2n+4)(2n+5), \text{ and}\\
\rho(\bl,-1_k) &= \rho(\bl) (2n+4-k) (2n+5-k) \cdots (2n+2) (2n+3) \\
&\text{ (just $\rho(\bl)$ if $k=0$)} 
\end{align*}
\item{(ii)} If $k\ge 1$ and $(\bl, 0,-1_k)\in\M_n$, then
$$
\rho(\bl,1,-1_k) = \rho(\bl) (2n+5-k) (2n+6-k) \cdots (2n+4)(2n+5).
$$
\end{prop}
\begin{proof} If $(\bl,0,-1_k)\in\M_{n+1}$ then $\rho(\bl,0,-1_k) = \rho(\bl)\alpha_1\alpha_2\cdots\alpha_{k+1}$, with
$\alpha_{k+1} = 2n+5$ and $\alpha_{j+1}-\alpha_j=1$ for each $j$, so
$$
\rho(\bl,0,-1_k) = \rho(\bl) (2n+5-k) (2n+6-k) \cdots (2n+4)(2n+5).
$$
The rest is similar.
\end{proof}

\bigskip

\begin{teor}\label{BerMot} For each $n\ge3$, let
$$
a_{ij}^{(n)} = \begin{cases}\dfrac{(j+1)(4n-j+1)}8 &\text{ if }j \text{ is odd and }i\le j-1\,,\\
              \dfrac{(j+2)(4n-j)}{32} &\text{ if }j \text{ is even and }i= j-1\,,\\
            0 &otherwise,\end{cases}
$$
and let $A^{(n)}=(a_{ij}^{(n)})_{1\le i \le n-2,\ 1\le j \le n-1}$.

Then 
$$
\zeta(2n) = \frac{\pi^{2n}}{(2n)!}\cdot\frac{4^{n-1}}{4^n-1}\, b_n\,,
$$
where $b_n$ is the first component of $A^{(3)}A^{(4)}\cdots A^{(n)}$. 

Therefore
\begin{align*}
B_{2n} &= (-1)^{n+1} \frac{b_n}{2(4^n-1)}\,,\text { and }\\
\tan^{(2n-1)}(0) &= \frac{4^{n-1}}n\,b_n\,.\\
\end{align*}
\end{teor}

\noindent{\bf Remark.} Note that the matrices $A^{(n)}$ are as follows:
$$
A^{(3)} = \begin{pmatrix} 3 & 5/4 \end{pmatrix},\qquad A^{(4)} = \begin{pmatrix} 4 & 7/4 & 0  \\ 4 & 0 & 7 \end{pmatrix},
$$
$$
 A^{(5)} = \begin{pmatrix} 5 & 9/4 & 0 & 0 \\ 5 & 0 & 9 & 0 \\ 5 & 0 & 9 & 12/4 \end{pmatrix},\qquad
A^{(6)} = \begin{pmatrix} 6 & 11/4 & 0 & 0 & 0 \\ 6 & 0 & 11 & 0 & 0 \\ 6 & 0 & 11 & 15/4 & 0 \\ 6 & 0 & 11 & 0 & 15
\end{pmatrix},  
$$
$$
A^{(7)}= \begin{pmatrix} 7 & 13/4 & 0 & 0 & 0 & 0 \\ 7 & 0 & 13 & 0 & 0 & 0 \\ 7 & 0 & 13 & 18/4 & 0 & 0 
\\ 7 & 0 & 13 & 0 & 18 & 0 \\  7 & 0 & 13 & 0 & 18 & 22/4 \end{pmatrix}
$$
and so on.

\begin{proof} Write formula (6) as
$$
\zeta(2n) = \frac {\pi^{2n}}{6(4^n-1)}\,\sum_{\bl\in\M_{n-2}} f(\bl)\,,
$$
with $f(\bl) = \frac{\nu(\bl)}{\rho(\bl)}$ for each $\bl \in \M$.

From Propositions \ref{general} and \ref{ro} it follows that
\begin{align*}
f(\bl,0,-1_k) &= \frac {2 (k+1) (2n+4-k)}{(2n+4) (2n+5)}\, f(\bl,-1_k)  \text { if $(\bl,-1_k)\in\M_n$, and}\\
f(\bl,1,-1_k) &= \frac {(k+1) (2n+4-k)}{2(2n+4) (2n+5)}\, f(\bl,0,-1_{k-1})  \text { if $(\bl,0,-1_{k-1})\in\M_n$.}
\end{align*}

Hence $f$ satisfies conditions $I$ and $II$ as in Theorem \ref{Mabs}, with 
$$b_{n,k} = \frac {2 (k+1) (2n+4-k)}{(2n+4) (2n+5)}\text{ and }c_{n,k}= \frac 1 4\, b_{n,k}\,.$$

Note that $f(0) = \dfrac 2 5$, and then Theorem \ref{Mabs} gives, for any $n\ge 3$,
$$
\frac 5 {2n+1} \sum_{\bl\in\M_{n-2}} f(\bl) = \langle \tilde A^{(3)} \tilde A^{(4)} \cdots \tilde A^{(n)}, e_1\rangle
$$
where $\tilde A^{(k+2)}$ is the $k$-th matrix in the statement of Theorem \ref{Mabs} and $\langle\cdot,\cdot\rangle$ is the
usual scalar product.

Let now $A^{(n)} = \dfrac {2n (2n+1)} 4 \tilde A^{(n)}$. It is easily checked that $A^{(n)}$ are the matrices in the statement,
and we have now
$$
\sum_{\bl\in\M_{n-2}} f(\bl) = \frac {6\cdot 4^{n-1}} {(2n)!}\, \langle  A^{(3)}  A^{(4)} \cdots  A^{(n)},
e_1\rangle\,,
$$
so
$$
\zeta(2n) = \pi^{2n} \,\frac {4^{n-1}}{4^n-1}\,\frac 1 {(2n)!} \,\langle  A^{(3)}  A^{(4)} \cdots  A^{(n)},
e_1\rangle\,.
$$
\end{proof}

\bigskip\bigskip\bigskip\bigskip\bigskip

\section{Entringer's theorem and Catalan numbers.}

\noindent In order to prove Theorem \ref{BerCat} we will make use of Entringer's theorem about alternating (or {\sl zig-zag\/})
permutations.

A permutation $\sigma\in S_n$ is said {\sl alternating\/} if it is such that
$$
\sigma(j) < \sigma(j+1) \text { if and only if } \sigma(j+1) > \sigma (j+2)
$$
for $j=1,\dots, n-2$, i.e. if $\sigma(j)$ is not a number between $\sigma(j-1)$ and $\sigma(j+1)$ for $j=2,\dots,n-1$.

Let $\alpha_n$ the number of alternating permutations. Let $\tau\in S_n$ given by $\tau(j) = n+1-j$, i.e. 
$$
\tau\equiv n,n-1,n-2,\dots,2,1.
$$
Then $\sigma\mapsto \tau\circ\sigma$ defines a bijection between alternating permutations $\sigma$ such that
$\sigma(1)<\sigma(2)$ and those such that $\sigma(1)>\sigma(2)$, and thus the number of any of these is $\alpha_n / 2$. 

Let then $\beta_n = \alpha_n/2$, with $\beta_0=\beta_1=1$. Entringer (see \cite{Entr}) proved that these numbers give a
combinatorial interpretation of tangent and secant numbers, namely
\begin{align*}
\sec z &= \beta_0 + \beta_2 \frac {z^2}2! + \beta_4 \frac{z^4}{4!} +\cdots \quad \text{ and}\\
\tan z &= \beta_1 z + \beta_3 \frac {z^3}3! + \beta_5 \frac{z^5}{5!} +\cdots 
\end{align*}
for each $z \in \C$ with $|z|<\pi/2$. In particular
$$
\beta(2n-1) = \tan^{(2n-1)}(0)
$$
for all $n\in\N$.

Independently from Calabi's argument, in \cite{Stanley} R. Stanley  obtained the tangent part by considering the polytopes in
(1)
 and, as we mentioned in the introduction, N. D. Elkies (\cite{Elkies}) has derived the result for both secant and tangent
numbers starting from Calabi's idea.

\bigskip

As for tangent numbers, there is the following relation with $\Phi$:

\begin{prop} Given $n\in\N$ and $\sigma\in S_{2n+1}$, $\Phi(\sigma)$ is a Dyck path if and only if $\sigma$ is
alternating and $\sigma(1) > \sigma(2)$.
\end{prop}
\begin{proof} $\Phi(\sigma) \in \D_{2n,0}$ means that $\lambda_j=\pm1$ for all $j$, so in
$$
0,\sigma(1),\sigma(2),\sigma(3),\dots,\sigma(2n),\sigma(2n+1),0
$$ 
or two or none of the neighbours of each $j$ are $>j$. This is exactly as saying that $\sigma$ is alternating and
$\sigma(1)>\sigma(2)$ (the fact that $\sigma(2n)<\sigma(2n+1)$ follows from them). 
\end{proof}

\begin{coro}\label{alter} For each $n\in \N$
$$
\tan^{(2n+1)}(0) = \sum_{\bl \in \D_{2n,0}} \nu(\bl)\,.\eqno{\qed}
$$
\end{coro}

\bigskip\bigskip

\noindent{\bf Proof of Theorem 1.1} By Proposition \ref{general}, for each $k\ge 0$
$$
\nu(\bl,1,-1_{k+1}) = (k+1) (k+2) \nu(1,-1_k)
$$
whenever $(\bl,-1_k)\in\D_0$. We can use Theorem \ref{Cabs} with $f = \nu$, and this together with Corollary
\ref{alter} gives the result.\hfill\qed

\bigskip\bigskip

Let's see another application of Theorem \ref{Mabs}, giving yet another way to obtain tangent (Bernoulli)
numbers by means of triangles.

\begin{teor}\label{nuzeros} Let $T(x)$ be the triangle generated by the infinite matrix
$$
A(x) = \begin{pmatrix}
2x & 1/x & 0 & 0 & 0 & 0 & \cdots  \\
2x & 0 & 4x & 0 & 0 & 0 &  \cdots  \\
2x & 0 & 4x & 3/2x & 0 & 0 &  \cdots   \\
2x & 0 & 4x & 0 & 6x & 0 & \cdots  \\
2x & 0 & 4x & 0 & 6x & 2/x & \cdots   \\
\hdotsfor{7}
\end{pmatrix}\,
$$
(so that $A(1)$ is the matrix of Corollary \ref{factM}).
Then the first element of the $n$-th row of $T(x)$ is a polynomial $P_n(x)$ of degree $n$ such that $P_n(1) = (n+1)!$ and
$$
P_n(0) = \tan^{(n+1)}(0).
$$
\end{teor}
\begin{proof} Let $f\colon\M\to\R$ given by $f(\bl) = \nu(\bl) x^{z(\bl)}$, with $z(\bl)$ the number of 0's in $\bl$, as in
Theorem
\ref{fceros}. By Proposition \ref{general} we have
\begin{align*} 
f(\bl, 0, -1_k) &= 2(k+1) \,x \,f(\bl, -1_k), \ \text{ and}\\
f(\bl, 1, -1_k) &= \frac{k+1}{2 x} \,f(\bl, 0, -1_{k-1}),
\end{align*}
and by Theorem \ref{Mabs}  the first element in the $n$-th row of $T(x)$ is 
$$
\sum_{\bl\in\M_n} f(\bl) = \sum_{\bl\in\M_n} \nu(\bl) x^{z(\bl)} = \sum_{k=0}^n \,u_{n,k} \,x^k =: P_n(x)\,,
$$ 
where $u_{n,k} = \d\sum_{\bl\in\D_{n,k}} \nu(\bl)$.

If $n$ is odd then $P_n(0)= 0$ since $\D_{n,0}$ is empty.

If $n$ is even, we have seen in Corollary \ref{alter} that $P_n(0) = \tan^{(n+1)}(0)$. This holds for $n=0$ too, since 
$P_0 = 1$.

Finally, $P_n(1) = (n+1)!$ as in Corollary \ref{factM}.
\end{proof}

\bigskip

\noindent{\bf Remark.} Of course, the matrix $A(x)$ is not defined for $x=0$, but we certainly have that
$$
\lim_{x\to 0}P_n(x) = \tan^{(n+1)}(0).
$$
If $x$ is positive and small enough then $\tan^{(n+1)}(0)$ will equal the floor function applied to $P_n(x)$. If we fix $n$ and
want to find $\tan^{(n+1)}(0)$ it suffices to take $x = 1/ (n+1)!$: we have seen that
$$
P_n(x) - \tan^{(n+1)}(0) = \sum_{k=1}^n \big(\sum_{\bl\in\D_{n,k}} \nu(\bl)\big) x^k\,,
$$
so if $0<x<1$ and $n$ is even
\begin{align*}
0<P_n(x) - \tan^{(n+1)}(0) &< x \,\sum_{k=1}^n \big(\sum_{\bl\in\D_{n,k}} \nu(\bl)\big)\\
&< x \,P_n(1) = x\,(n+1)!\,,
\end{align*}
and for $x = \dfrac 1 {(n+1)!}$ we have $0<P_n(x) - \tan^{(n+1)}(0) < 1$.

\bigskip\bigskip

Entringer's paper \cite{Entr} shows actually more than we have said so far: it provides a triangle to generate both tangent
and secant numbers. Our last aim is to recall this, presenting it in the frame of our matrix-generated triangles idea.  

Let $E_{n,k}$ denote the number of alternating permutations $\sigma\in S_{n+1}$ such that $\sigma(1) = k+1>\sigma(2)$,
with $E_{0,0}=1$. Note that $E_{n,0} = 0$ for $n\ge 1$, and $E_{n,n}$ is just the zig-zag number $\beta_n$.

$(E_{n,k})$ are the so-called {\sl Entringer numbers\/}, and they form the {\sl Seidel-Entringer-Arnold triangle\/}. The
starting point for the proof of Entringer's theorem in
\cite{Entr} is the recurrence
$$
E_{n+1,k+1} = E_{n+1,k} + E_{n,n-k} \quad (n\ge k\ge 0)\,,
$$
which gives itself a very simple algorithm to generate the triangle and thus the zig-zag numbers (see \cite{Sloane} for more on
this). This recurrence yields
\begin{equation} \label{Ent}
E_{n+1,k+1}=\sum_{j=n-k}^n E_{n,j} \quad (n\ge k\ge 0)\,,
\end{equation}
a fact that can be formulated as the following observation: if we remove the first (and trivial) column in the
Seidel-Entringer-Arnold triangle, the resulting one is easily obtained by means of matrices:

\begin{prop} For each $n\in\N$, let $(A^{(n)} = (a_{ij})$ the $n\times (n+1)$ matrix given by
$$
a_{ij} = \begin{cases} 0 &\text{ if }i+j\le n\,,\\
              1 &\text{ if }i+j>n\,.\end{cases}
$$
Then the triangle $(t_{n,m})_{0\le m\le n)}$ generated by $(A^{(n)})$ is such that
$$
t_{n,m} = E_{n+1,m+1} \quad \text{ for all }n,m\,.
$$
In particular $t_{n,0}$ is the zig-zag number $\beta_n$ for all $n\ge 0$, and $t_{n,n} =\beta_{n+1}$.
\hfill\qed
\end{prop}

\bigskip

Note that $A^{(n)}$ in the proposition is the submatrix of the first $n$ rows and the {\sl last\/} $n+1$ columns of
the infinite matrix
$$
A = \begin{pmatrix}
\cdots & 0 & 0 & 0 & 1 & 1\\
 \cdots & 0 & 0 & 1 & 1 & 1 \\
\cdots & 0 & 1 & 1 & 1 & 1  \\
\hdotsfor{6}\\ \end{pmatrix}\,,
$$
and the triangle begins
$$
\begin{array}{cccccc}
1 &   &   &   &   &   \\
1 & 1 &   &   &   &   \\
1 & 2 & 2 &   &   &   \\
2 & 4 & 5 & 5 &   &   \\
5 & 10 & 14 & 16 & 16 &  \\
16 & 32 & 46 & 56 & 61  & 61 \\
\hdotsfor{6}
\end{array}
$$
\bigskip

\begin{coro} Let the infinite matrix
$$
A = \begin{pmatrix}
1 & 2 & 2 & 2 & 2 & 2 & 2 & \cdots  \\
1 & 2 & 3 & 3 & 3 & 3 & 3 &  \cdots  \\
1 & 2 & 3 & 4 & 4 & 4 & 4 &  \cdots  \\
1 & 2 & 3 & 4 & 5 & 5 & 5 &  \cdots  \\
1 & 2 & 3 & 4 & 5 & 6 & 6 & \cdots  \\
1 & 2 & 3 & 4 & 5 & 6 & 7 & \cdots  \\
\hdotsfor{8} \end{pmatrix}\,,
$$
and let $A^{(n)}$ be the submatrix formed by the first $n$ rows and the first $n+2$ columns of $A$, for all $n\in\N$. Let then 
$\bar t_1 = (1)$ and 
$$\bar t_{n+1} = (t_{n+1,1}, t_{n+1,2},\dots,t_{n+1,2n+1}) = \bar t_{n} A^{(2n-1)} \in \N^{2n+1} \quad (n\in\N)\,.
$$ Then, for every $n = 1,2,\dots$ we have
\begin{align*}
t_{n,1}&= \beta_{2n-2}\quad \text{ (the $n$-th secant number), and}\\
t_{n,2n-1} &= \beta_{2n-1} \quad \text{ (the $n$-th tangent number).}
\end{align*}
\end{coro}
\begin{proof} Just note that, if $\tilde A^{(n)}$ are the matrices in the previous proposition, then 
$$
\tilde A^{(2n-1)} \tilde A^{(2n)} = A^{(2n-1)}
$$ 
for each $n\in\N$, and thus the vectors $\bar t_n$ are the even rows in the triangle of the previous proposition.
\end{proof}

\begin{prop} Let $A=(a_{ij})$ the infinite matrix given by
$$
a_{ij} = \begin{cases} 1 &\text{ if }i \text{ is odd and }i\le j\,,\\
              1 &\text{ if }j\text{ is odd and }j\le i+1\,,\\
               0 &\text{ otherwise,}
\end{cases}
$$
i.e.
$$
A = \begin{pmatrix}
1 & 1 & 1 & 1 & 1 & 1 & 1 & \cdots  \\
1 & 0 & 1 & 0 & 0 & 0 & 0 &  \cdots  \\
1 & 0 & 1 & 1 & 1 & 1 & 1 &  \cdots  \\
1 & 0 & 1 & 0 & 1 & 0 & 0 &  \cdots  \\
1 & 0 & 1 & 0 & 1 & 1 & 1 & \cdots  \\
1 & 0 & 1 & 0 & 1 & 0 & 1 & \cdots  \\
\hdotsfor{8} \end{pmatrix}\,.
$$
Then the triangle $(t_{n,m})_{n\ge m \ge0}$ generated by $A$ is such that
$$
t_{n,0} = \beta_{n+1} \quad \text{ (the zig-zag number)}
$$ 
and each row $\bar t_n$ is a permutation of the Entringer numbers $(E_{n+1,k})_{k=1,\dots,n+1}$.
\end{prop}
\begin{proof} Using (7), it follows (by induction on $n$) that
\begin{align*}
t_{n,2k} &= E_{n+1, n+1-k} \quad \text{ if $0\le k \le \frac n 2$, and}\\
t_{n,2k+1}&= E_{n+1, k+1} \quad \text{ if $0\le k \le \frac{n-1}2$,}
\end{align*}
and then
$$
(t_{n,m})\equiv 
\begin{array}{cccccc}
1 &   &   &   &   &   \\
1 & 1 &   &   &   &   \\
2 & 1 & 2 &   &   &   \\
5 & 2 & 5 & 4 &   &   \\
16 & 5 & 16 & 10 & 14 &  \\
61 & 16 & 61 & 32 & 56  & 46 \\
\hdotsfor{6}
\end{array}
$$
\end{proof}

\bigskip\bigskip



\begin{thebibliography}{99}

\bibitem{Aig} M. Aigner, 
\textit{Motzkin numbers},
{Europ. J. Combinatorics, 19}
 (1998), 663--675.

\bibitem{Barc} E. Barcucci, R. Pinzani and R. Sprugnoli,
\textit{The Motzkin family},
{Pure Math. and App., Ser. A, 2} (1991), 249--279.

\bibitem{Calabi} F. Beukers, J. A. C. Kolk and E. Calabi,
\textit{Sums of generalized harmonic series and volumes},
{Nieuw Arch. Wisk. (4), 11} (1993), 217--224.


\bibitem{Dyck} E. Deutsch,
\textit{Dyck path enumeration},
{Discrete Math., 204} (1999), 167--224.

\bibitem{Shap} R. Donaghey and L. W. Shapiro,
\textit{Motzkin numbers},
{J. Combin. Th. Ser. A, 23} (1977), 291--202.

\bibitem{Elkies} N. D. Elkies,
\textit{On the sums $\d\sum_{k=-\infty}^{\infty}(4k+1)^{-n}$},
{preprint, avalaible in the World Wide Web at arXiv:math.CA/0101168 v2}(2001).

\bibitem{Entr} R. C. Entringer,
\textit{A combinatorial interpretation of the Euler and Bernoulli numbers},
{Nieuw Arch. Wisk. (3), 14} (1966), 241--246.

\bibitem{Gardner} M. Gardner,
\textit{Catalan numbers: An integer sequence that materializes in unexpected places},
{Sci. Amer., 234}(1976), 120--125.

\bibitem{GKP} R. L. Graham, D. E. Knuth and O. Patashnik,
\textsl{Concrete Mathematics},
Addison-Wesley, 1994.

\bibitem{Hardy} G. H. Hardy and E. M. Wright,
\textsl{An Introduction to the Theory of Numbers},
Oxford Science Publications, 1979.

\bibitem{Sloane} J. Millar, N. J. A. Sloane  and N. E. Young,
\textit{A new operation on sequences: the boustrophedon transform},
{J. Combin. Th. Ser. A, 76} (1996), 44--54.


\bibitem{Stanley} R. Stanley, \textit{Two poset polytopes}, {Discrete and Computational Geometry, 1} (1986), 9--23. 




\end{thebibliography}
\enddocument